\documentclass[11pt]{article}

\usepackage{graphicx}

\newfont{\bbb}{msbm10 scaled\magstephalf}
\newfont{\sbbb}{msbm7 scaled\magstephalf}
\newfont{\frak}{eufm10 scaled\magstep1}
\def\C{\mbox{\bbb{C}}}
\def\P{\mbox{\bbb{P}}}
\def\Q{\mbox{\bbb{Q}}}
\def\R{\mbox{\bbb{R}}}
\def\Z{\mbox{\bbb{Z}}}

\def\e{e^{2\pi i\theta}}
\def\eh{e^{2\pi i(\theta+h)}}
\def\ep{e^{2\pi ip\theta}}
\def\eq{e^{2\pi iq\theta}}
\def\es{e^{2\pi is\theta}}
\def\et{e^{2\pi it\theta}}
\def\esh{e^{2\pi is(\theta+h)}}
\def\eth{e^{2\pi it(\theta+h)}}
\def\zs{|z|^2}
\def\ws{|w|^2}

\title{Toric Quasifolds}
\author{Elisa Prato\thanks{Partially supported by the PRIN Project ``Real and Complex Manifolds: Topology, Geometry and Holomorphic Dynamics" (MIUR, Italy) and by GNSAGA (INdAM, Italy).}}

\date{}
\begin{document}
\maketitle
\section*{}

Quasifolds are a class of highly singular spaces. They are locally modeled 
by manifolds modulo the smooth action of countable groups. If the countable groups 
happen to be all finite, then quasifolds are orbifolds and if they happen to be all equal to the identity, they are manifolds. 
They were first introduced in \cite{p} in order to address, from the symplectic viewpoint,
the longstanding open problem of extending the classical constructions of toric geometry
to those convex polytopes that are \textit{not} rational. 

In order to clarify this last statement, let us begin by recalling what it means for a convex polytope to be rational. 
It is well known that every convex polytope in $(\R^n)^*$  can be written as the bounded
intersection of finitely many closed half--spaces:
\begin{equation}
\label{decomp}\Delta=\bigcap_{j=1}^d\{\;\mu\in(\R^n)^*\;|\;\langle\mu,X_j\rangle\geq\lambda_j\;\},
\end{equation}
where $X_1, \ldots X_d\in \R^n$, $\lambda_1,\ldots,\lambda_d\in \R$, and $d$ is the number of
\textit{facets} (codimension--one faces) of $\Delta$
\cite[Theorem~1.1]{ziegler}. It is not restrictive to assume that $\Delta$ has full dimension $n$.
We remark that the vectors $X_1, \ldots X_d$
are orthogonal to the facets of $\Delta$ and inward--pointing. For brevity, we will be referring to 
these vectors as \textit{normals} for $\Delta$.
The polytope is then said to be \textit{rational} if the normals can be chosen inside of a
lattice $L\subset\R^n$. Rationality is a rather restrictive condition, and, in fact, many interesting
convex polytopes are not rational: take, for instance, the regular pentagon and the regular dodecahedron.

Now, toric geometry, initiated by Demazure in
\cite{demazure}, sets to associate with each rational convex polytope
a beautiful geometrical space with special torus symmetries. One of the remarkable consequences of doing so is that
the geometry of the space can be used to deduce combinatorial information on the polytope and viceversa.
The construction of toric spaces can be done from different geometric perspectives: 
algebraic \cite{fulton}, complex \cite{audin,cox} and symplectic \cite{delzant}
\footnote{The starting point in the algebraic and complex category is actually, more generally, a {\em fan} instead of a polytope, 
but the basic idea that follows applies verbatim (see \cite{whatis}).}.
The crucial fact to recall here is that these constructions
always rely on the lattice $L$ and on a set of primitive normals 
in $L$. Evidently, for nonrational polytopes this setup is missing.
The first step in generalizing toric geometry to this case (see \cite{p}) consists in replacing the lattice 
with a similar enough object, which allows sufficient freedom to contain 
a set of normals for the polytope. The optimal choice turns out to be that of a
\textit{quasilattice} $Q$, namely the $\Z$--span of a set of $\R$--spanning vectors of $\R^n$.  
\begin{figure}[h]\begin{center}
\includegraphics[scale=0.6]{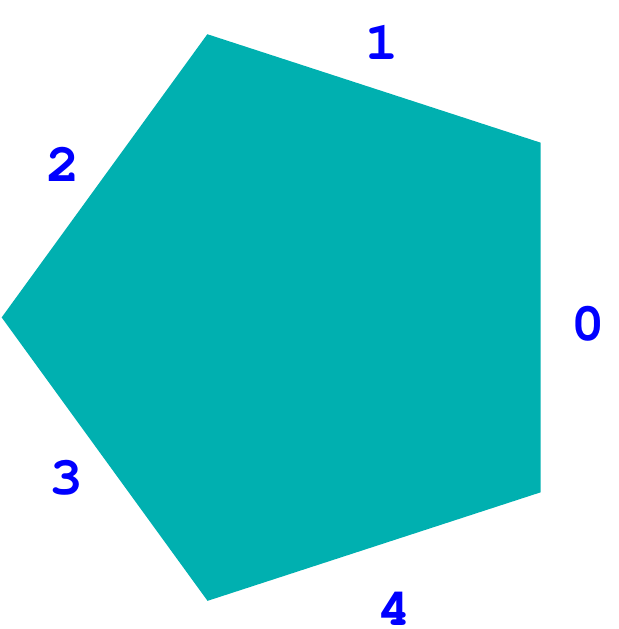} \quad\quad
\includegraphics[scale=0.3]{star}
\caption{The regular pentagon and the fifth roots of unity}
\label{pentagonandstar}
\end{center}
\end{figure}
In the case of the regular pentagon, for example, one considers the $\Z$--span of the 
fifth roots of unity (see Figure~\ref{pentagonandstar}). 
We thus have a new framework given by the triple
\begin{equation}
\label{triple}
(\Delta, Q, \mbox{normals in } Q)
\end{equation}
and, once this has been fixed, the standard toric constructions can be extended. 
For polytopes that are \textit{simple} (meaning that each vertex is the intersection of exactly $n$ facets),  they give rise to what we call \textit{toric quasifolds}. This was done first in the symplectic category 
\cite{p} and then, jointly with Battaglia, in the complex/K\"ahler category \cite{cx}.
The torus symmetries of the rational case are replaced by the symmetries of a {\em quasitorus}: it is the abelian
group $\R^n/Q$, which is itself a quasifold.
Though not Hausdorff in general, toric quasifolds have beautiful atlases that generalize the 
standard toric atlases of the rational case: each chart is the quotient of an 
open subset of $\C^n\simeq \R^{2n}$ modulo the smooth action of a \textit{countable} subgroup of the standard torus $\R^n/\Z^n$.

Battaglia has extended both the symplectic and complex/K\"ahler constructions to completely general convex polytopes, 
no longer necessarily simple; the resulting toric spaces are even more singular, but they turn out to be
stratified by toric quasifolds \cite{stratif-re, stratif-cx}.

It is interesting to remark that quasilattices are also crucial in the theories of nonperiodic tilings
(see \cite{mackay1} and \cite[Chapter 2]{senechal1}). The pentagonal quasilattice above, for example, 
arises in relation with Penrose tilings.

It is our goal here to illustrate toric quasifolds, and their atlases, by describing a number of examples. 
We do so in the symplectic category, but of course everything can be reformulated in the complex one.
We begin with a $2$--dimensional example that displays all of the main characteristics of quasifolds: the {\em quasisphere}. 
We pass on to considering examples of dimension $4$ and $6$
that came about by exploring the natural connection with Penrose and Ammann tilings. 
We then briefly address the toric spaces corresponding to
the regular convex polyhedra. We conclude with a number of considerations.

For the formal definition of quasifold, we refer the reader to \cite{p,cx}. 
The complex and symplectic atlases for toric quasifolds are explicitly described in \cite[proofs of Theorems~2.2 and~3.2]{cx}.

\subsection*{From sphere to orbisphere to quasisphere}
Quasispheres, introduced in \cite{p}, are generalizations of spheres and orbispheres, 
so we will begin by recalling some relevant facts on the latter two.

For any positive real number $r$, let $B(r)\subset \C$ be the open ball of 
center the origin and radius $\sqrt{r}$. Consider, for any positive real number $\alpha$, the group
$$\Gamma_{\alpha}=\{\, e^{2\pi ik\alpha}\in S^1\,|\, k\in \Z\,\}.$$
Notice that $\Gamma_{\alpha}$ is the identity when $\alpha$ is an integer, it is finite for $\alpha$ rational, while it is
countable for $\alpha$ irrational. The group $\Gamma_{\alpha}$ acts on the ball $B(r)$ by complex multiplication. 
For any $z\in B(r)$, we will denote by $[z]\in B(r)/\Gamma_{\alpha}$ the corresponding orbit.

\subsubsection*{The sphere}
Let us write the $2$ and $3$--dimensional unit spheres as follows
$$S^2=\{\,(z,x)\in \C\times\R\,|\, \zs+x^2=1\,\},$$
$$S^3=\{\,(z,w)\in \C^2\,|\, \zs+\ws=1\,\}.$$
The surjective mapping
\begin{eqnarray*}
f\,\colon &S^3 &\longrightarrow S^2\\
&\left(z,w\right)&\longmapsto \left(2z\overline{w}, \zs-\ws\right)
\end{eqnarray*}
is known as the {\em Hopf fibration}. It is easily seen that the fibers of this mapping are given by the orbits of the 
circle group
$$S^1=\{\,\e\,|\,\theta\in\R\,\}$$
acting on $S^3$ by complex multiplication as follows: 
$$\e\cdot(z,w) = \left( \e z, \e w\right).$$
Therefore $S^2$ can be identified with the space of orbits $S^3/S^1$. 
Notice that the $S^1$--orbits through the points $(0,1)$ and $(1,0)$ of $S^3$ correspond, respectively, 
to the south pole, $S=(0,-1)$, and north pole, $N=(0,1)$, of $S^2$.

For each $(z,w)\in S^3$, we denote by $[z:w]\in S^3/S^1\simeq S^2$ the corresponding orbit.
Let us describe the standard atlas of $S^2$. Consider the covering given by the open subsets 
$$U_S=\{\,[z:w]\in S^2\;|\;w\neq 0\,\}$$
$$U_N=\{\,[z:w]\in S^2\;|\;z\neq 0\,\}.$$ As the notation suggests, the first is a neighborhood of the south pole $S=[0:1]$, 
while the second is a neighborhood of the north pole $N=[1:0]$. Finally, we have homeomorphisms:
$$
\begin{array}{ccc}
B(1) &\longrightarrow &U_S\\
z&\longmapsto &\left[z:\sqrt{1-\zs}\right]
\end{array}
$$
$$
\begin{array}{ccc}
B(1) &\longrightarrow &U_N\\
w&\longmapsto &\left[\sqrt{1-\ws}:w\right].
\end{array}
$$

\subsubsection*{The orbisphere}
This simple quotient construction can be extended to the orbifold setting as follows. Let $p,q$ be two relatively 
prime positive integers and consider the $3$--dimensional ellipsoid $$S^3_{p,q}=\{\,(z,w)\in \C^2\,|\, p\zs+q\ws=pq\,\}.$$
The circle group $S^1$ acts on $S^3_{p,q}$ as follows:
\begin{equation}\label{circleaction}\e\cdot (z,w) = \left( \ep z, \eq w\right).\end{equation} 
Taking the space of orbits in this case yields the $2$--dimensional orbifold 
$S^2_{p,q}=S^3_{p,q}/S^1$, called {\em orbisphere}. It admits the two singular points $S=[0:\sqrt{p}]$ and $N=[\sqrt{q}:0]$.

Similarly to what we have done for the sphere, for each $(z,w)\in S^3_{p,q}$,
we denote by $[z:w]\in S^2_{p,q}$ the corresponding orbit. We then
consider the covering given by the two open subsets
$$U_S=\{\,[z:w]\in S^2_{p,q}\;|\;w\neq 0\,\}$$
$$U_N=\{\,[z:w]\in S^2_{p,q}\;|\;z\neq 0\,\}.$$
The first is a neighborhood of the point $S=[0:\sqrt{p}]$, while the second is a neighborhood of the point $N=[\sqrt{q}:0]$.
The mappings
$$
\begin{array}{ccc}
B(q)/\Gamma_\frac{1}{q}&\longrightarrow &U_S\\
\left[z\right]&\longmapsto &\left[z:\sqrt{p-\frac{p}{q}\zs}\right];
\end{array}
$$
$$
\begin{array}{ccc}
B(p)/\Gamma_\frac{1}{p}&\longrightarrow &U_N\\
\left[w\right]&\longmapsto &\left[\sqrt{q-\frac{q}{p}\ws}:w\right]
\end{array}
$$
are homeomorphisms, turning $U_S$ and $U_N$ into orbifold charts.
\subsubsection*{The quasisphere}
We now extend the construction even further. Let $s,t$ be two positive real numbers with $s/t \notin \Q$ and
consider the $3$--dimensional ellipsoid $$S^3_{s,t}=\{\,(z,w)\in \C^2\,|\, s\zs+t\ws=st\,\}.$$ 
Simply substituting $p,q$ with $s,t$ in (\ref{circleaction}), 
does not define an $S^1$--action on $S^3_{s,t}$: 
in fact, if you replace $\theta$ by $\theta +h$, where $h$ is a non--zero integer, we have  $\eh=\e$ but
$ (\esh,\eth)\neq(\es,\et)$.
The idea is to consider the irrational wrap on the standard two--torus instead:
$$N=\{\, (\es,\et)\in \R^2/\Z^2\,|\, \theta\in \R\,\}.$$
The standard action of $N$ on $S^3_{s,t}$ is now well defined and we take our {\em quasisphere} 
to be the space of orbits $S^2_{s,t}=S^3_{s,t}/N$. 
This quotient is the simplest example of quasifold. It is wilder then the sphere and orbisphere, 
in that it is not a Hausdorff topological space. However, quasisphere charts
are a straightforward and very natural generalization of the standard sphere and orbisphere charts. 

Exactly as done above, for each $(z,w)\in S^3_{s,t}$,
we denote by $[z:w]$ the corresponding orbit. We then
consider the covering of $S^2_{s,t}$ given by the opens subsets $$U_S=\{\,[z:w]\in S^3_{s,t}/\R\;|\;w\neq 0\,\}$$
$$U_N=\{\,[z:w]\in S^3_{s,t}/\R\;|\;z\neq 0\,\}.$$
The first is a neighborhood of the point $S=[0:\sqrt{s}]$, while the second is a neighborhood of the point $N=[\sqrt{t}:0]$.
They are each homeomorphic to the quotient of an open subset 
of $\C$ modulo the action of a {\em countable} group. 
In fact, the mappings
$$
\begin{array}{ccc}
B(t)/\Gamma_{\frac{s}{t}} &\longrightarrow &U_S\\
\left[z\right]&\longmapsto &\left[z:\sqrt{s-\frac{s}{t}\zs}\right]
\end{array}
$$
$$
\begin{array}{ccc}
B(s)/\Gamma_{\frac{t}{s}} &\longrightarrow &U_N\\
\left[w\right]&\longmapsto &\left[\sqrt{t-\frac{t}{s}\ws}:w\right]
\end{array}
$$
are homeomorphims. 

\noindent {\sc Remark.}
The sphere is the symplectic toric manifold corresponding to the unit interval, with lattice $\Z$ 
and primitive normals $X_1=1$, $X_2=-1$.  The orbisphere, on the other hand, is the symplectic toric orbifold corresponding
to the same interval, with same lattice and normals $X_1=q$, $X_2=-p$. Finally, the quasisphere is the symplectic toric quasifold 
corresponding to the same interval, with quasilattice $Q=s \Z+t\Z$ and normals $X_1=t$, $X_2=-s$.
Wanting to consider a rational polytope, such as the unit interval, in a nonrational setup may seem strange at first sight, but
in fact it is quite useful. We will see other instances of this in the next section. 
Also, the sphere and orbisphere provide the simplest examples showing that the same polytope and (quasi)lattice yield 
different symplectic toric spaces, if the normals are changed. The choice of normals within a same quasilattice is in fact totally free, 
but sometimes a natural choice is dictated by the context. This is actually the case for all of the examples that follow.

\subsection*{Quasifolds and nonperiodic tilings}
\subsubsection*{Quasifolds corresponding to Penrose and Ammann tilings}
The fact that quasilattices appear naturally in nonperiodic tilings
lead us to explore, jointly with Battaglia, the connection between toric quasifolds and Penrose
and Amman tilings.
 
\begin{figure}[h]
\begin{center}
\includegraphics[scale=0.3]{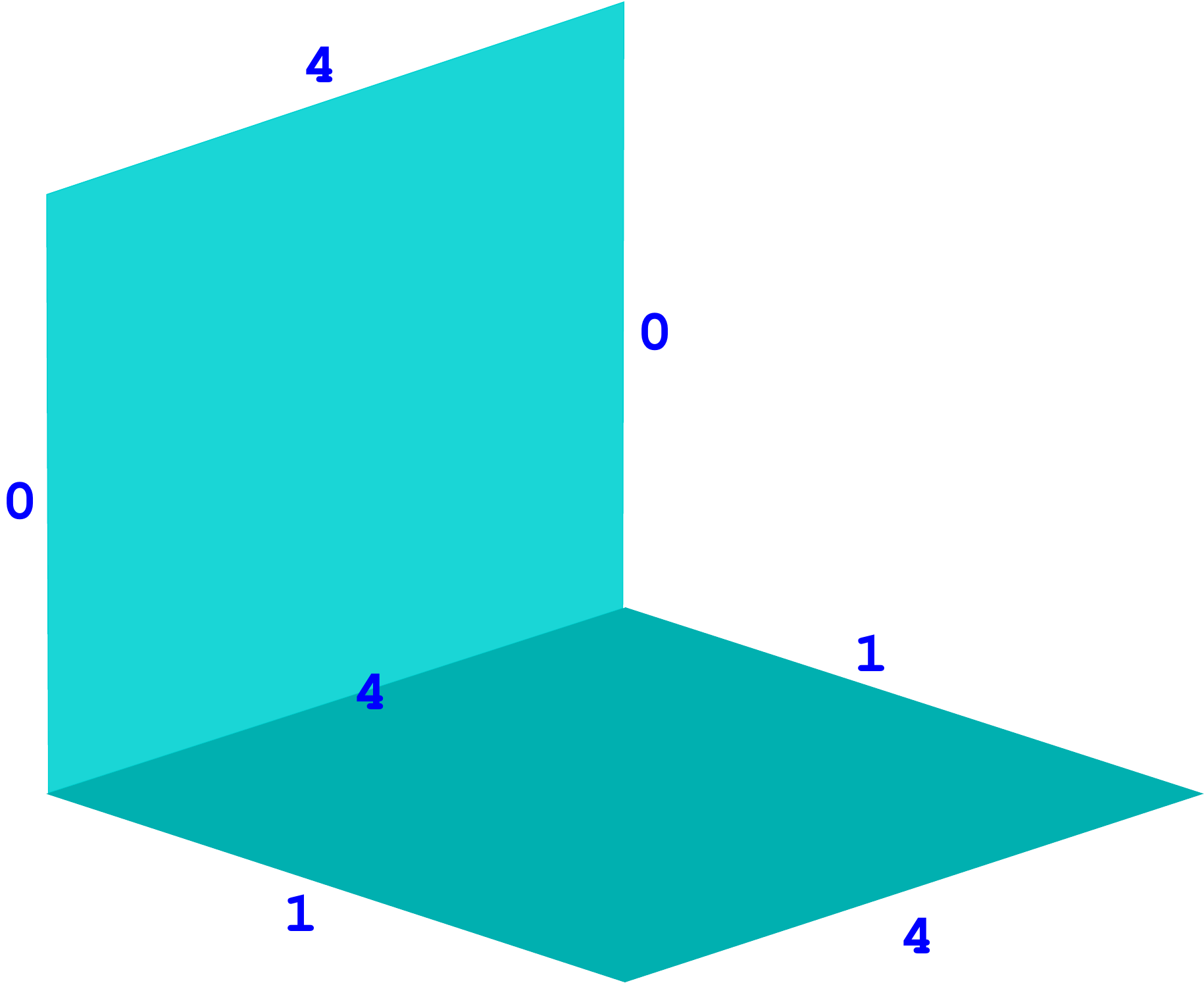} \quad\quad
\includegraphics[scale=0.3]{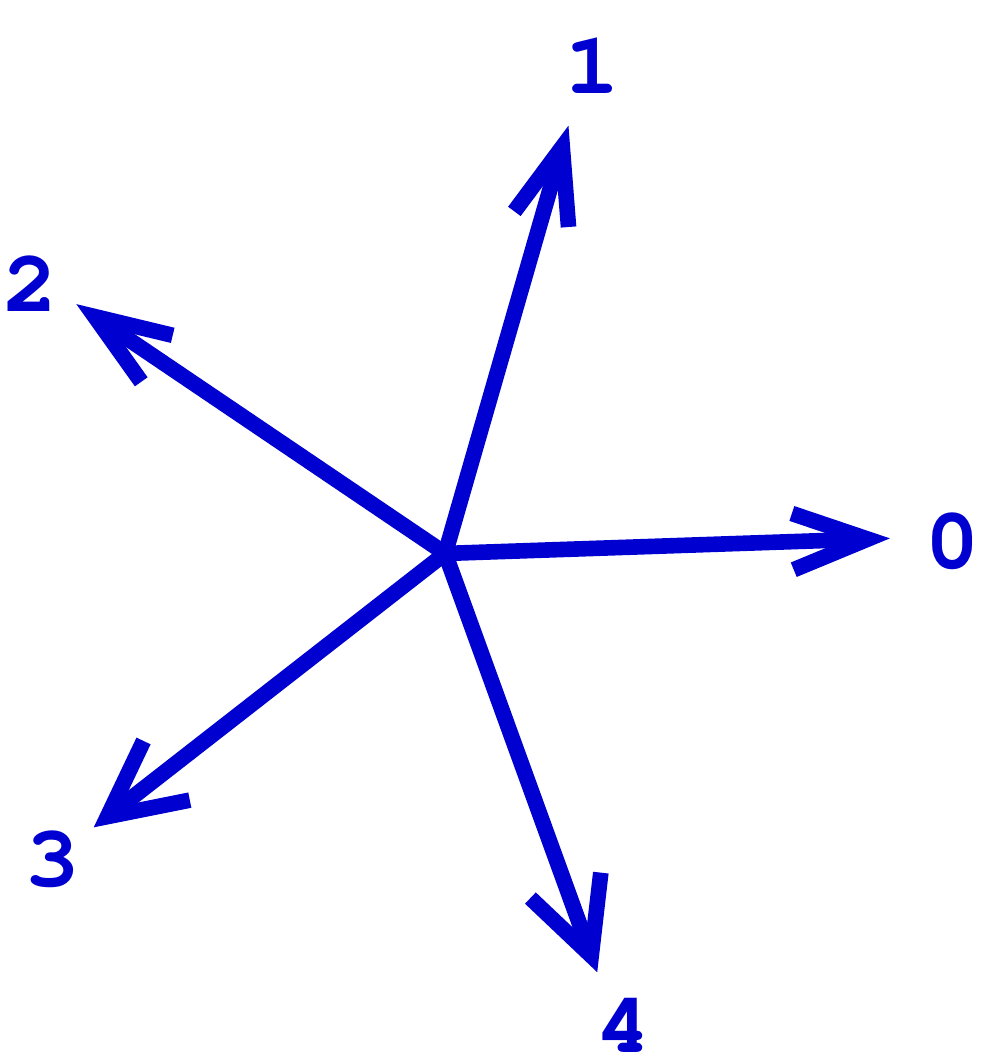}
\caption{The Penrose rhombuses}
\label{rhombusesandstar}
\end{center}
\end{figure}

 Penrose rhombus tilings are nonperiodic tilings that are composed by two 
 different types of rhombuses, thick and thin \cite{pentaplexity}.
 These rhombuses are simple convex polytopes and it is natural to want to compute the corresponding toric quasifolds.
 The normals of each rhombus taken separately actually span a lattice, so each of them is rational in its own right. 
 However, if we want to treat simultaneously all of the rhombuses of a given tiling, we need to consider a quasilattice: 
 the natural choice here is the pentagonal quasilattice that we introduced earlier,
with normals the relevant fifth roots of unity (see Figure~\ref{rhombusesandstar}).
The generalized toric construction then yields a pair of four--dimensional toric quasifolds, one for each different type of rhombus.
They are both given by a quotient of the type $(S^2(r)\times S^2(r))/\Gamma$, where $S^2(r)$ denotes the 
$2$--sphere of radius $r$ and $\Gamma$ is a countable subgroup of the standard 2--torus. 
The radius $r$ is  $(\frac{1}{2}\sqrt{2+\phi})^{1/2}$ for the thick rhombus and $(\frac{1}{2\phi}\sqrt{2+\phi})^{1/2}$ for the thin one,
where $\phi=\frac{1+\sqrt{5}}{2}$ is the \textit{golden ratio}. The two quasifolds are diffeomorphic but not symplectomorphic.

Something analogous happens for the three--dimensional generalization 
of this tiling due to Ammann, which is composed by two different types of rhombohedrons, prolate and oblate \cite{senechal2}.
Again, each rhombohedron is rational, but to treat them all of them simultaneously we need to
consider a quasilattice, known in crystallography as the {\em face--centered icosahedral lattice}. 
As normals we choose the relevant generators. 
One then obtains a pair of six--dimensional symplectic toric quasifolds, one for each 
type of rhombohedron. 
Similarly to what happens for the rhombus tiling, they are given by
$(S^2(r)\times S^2(r)\times S^2(r))/\Gamma$, where $\Gamma$ is a countable subgroup of the standard 3--torus.
The radius $r$ here is $[2\phi^2(3-\phi)]^{-\frac{1}{4}}$ for the oblate rhombohedron and $[2(3-\phi)]^{-\frac{1}{4}}$ 
 for the prolate one. Again, the two spaces here are diffeomorphic but not symplectomorphic.

As we have seen, the quasifolds for both Penrose rhombus tilings and Ammann tilings are \textit{global}, 
namely the quotient of a manifold modulo the action of a countable group.

\begin{figure}[h]
 \begin{center}
\includegraphics[scale=0.3]{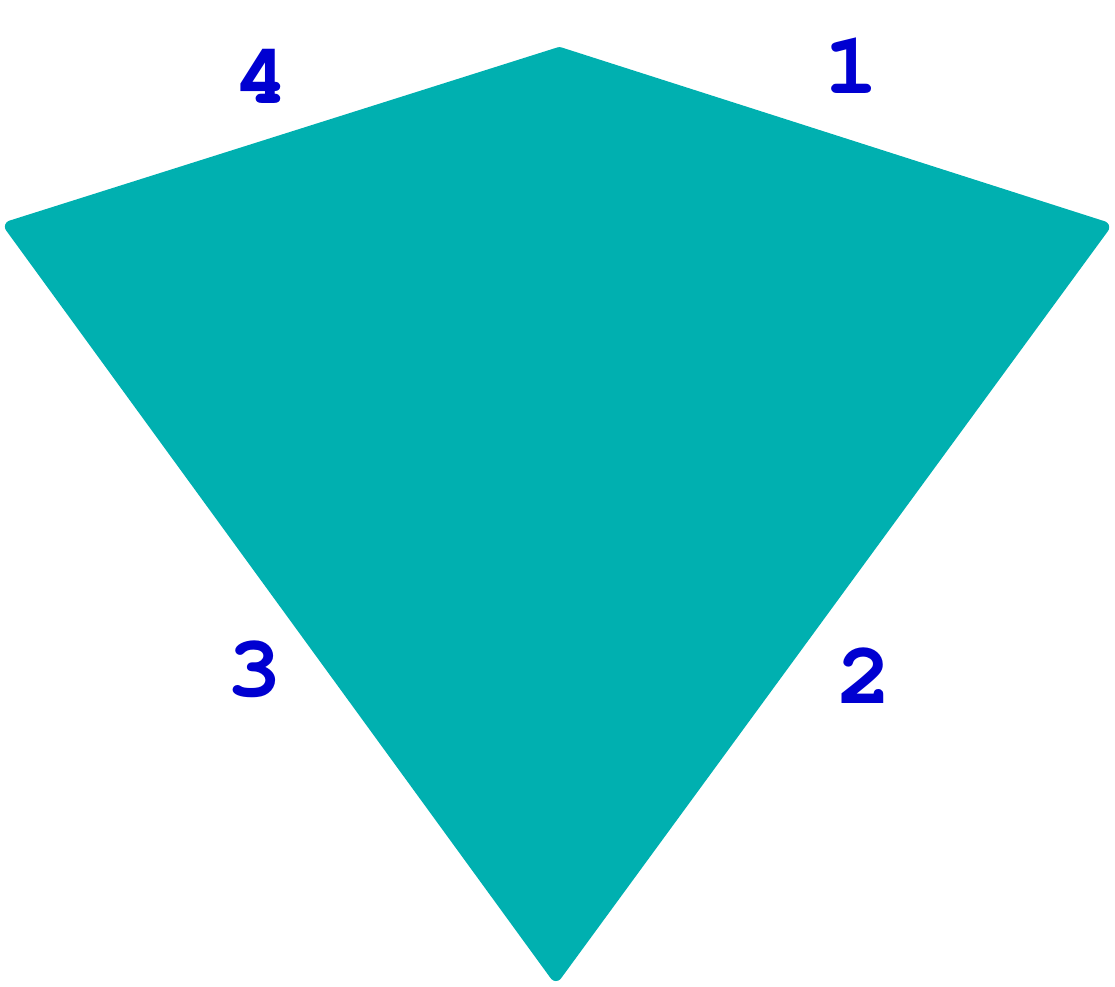} \quad\quad
\includegraphics[scale=0.3]{star.pdf}
\caption{The Penrose kite}
\label{kiteandstar}
\end{center}
\end{figure}

Something entirely different happens for the kite and dart tiling \cite{pentaplexity}. 
First of all, the only tile here that is convex, and therefore relevant
to our discussion, is the kite. Moreover, the kite, unlike the rhombuses and rhombohedrons, is actually nonrational. 
So there is no choice but to consider a quasilattice, and the natural one happens to be, again, the pentagonal quasilattice;
the normals are, up to sign, the relevant fifth roots of unity
(see Figure~\ref{kiteandstar}).
 Then the resulting toric quasifold is not global.
 It is the four--dimensional quasifold given by 
 $$
 M=\frac{{\left\{ \,(z_1,z_2,z_3,z_4)\in\C^4\,|\,\phi|z_1|^2+|z_2|^2+\phi |z_3|^2=\phi |z_2|^2+|z_3|^2+\phi |z_4|^2
=\phi\,\right\}}}{{\left\{\,\exp\left(-s+\phi t,s,t,-t+\phi s\right)\in \R^4/\Z^4\,|\,s,t \in\R\,\right\}}}.
 $$
Let us describe one of its charts.
Consider the open subset
$$\tilde{U}=\left\{\,(z_2,z_3)\in\C^2\,|\,|z_2|^2+\phi |z_3|^2<\phi,\,\phi |z_2|^2+|z_3|^2<\phi\,\right\}$$
and the countable group 
$$
\Gamma=\left\{\,(e^{2\pi i\phi h},e^{2\pi i\phi k})\in \R^2/\Z^2\;|\; h,k\in\Z\right\}.
$$
Then the mapping
$$
\begin{array}{ccc}
\tilde{U}/\Gamma & \longrightarrow & \left\{ [z_1:z_2:z_3:z_4]\in M\;|\;z_1\neq0,z_4\neq0\right\} \\
\left[z_2:z_3\right] & \longmapsto & \left[\sqrt{\phi-|z_2|^2-\phi |z_3|^2}:z_2:z_3:\sqrt{\phi-\phi |z_2|^2-|z_3|^2}\right]
\end{array}
$$
is a homeomorphism.

\subsubsection*{Decomposing Penrose tiles and symplectic cutting}
\begin{figure}[h]
\begin{center}
\includegraphics[scale=0.22]{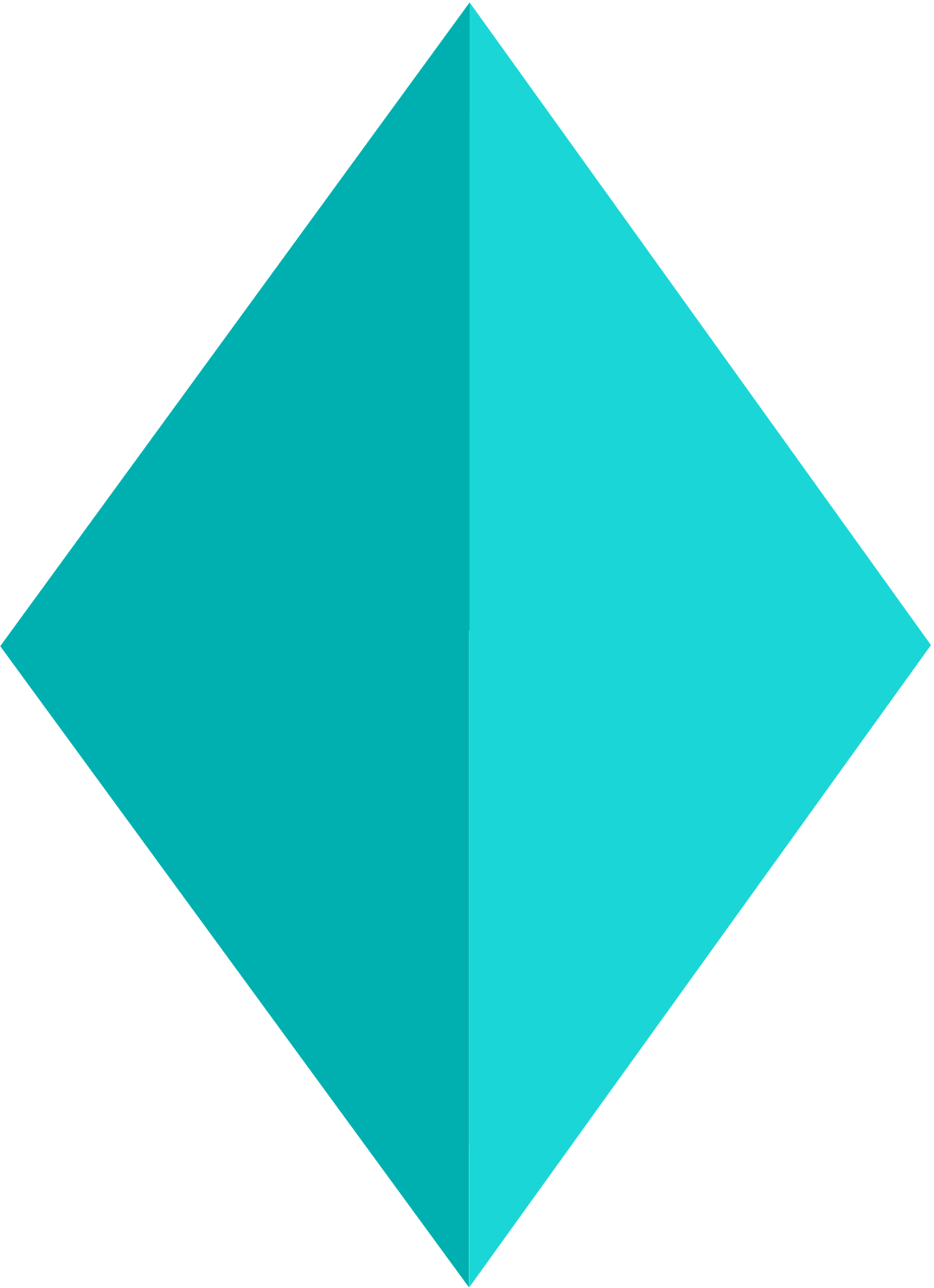}
\includegraphics[scale=0.22]{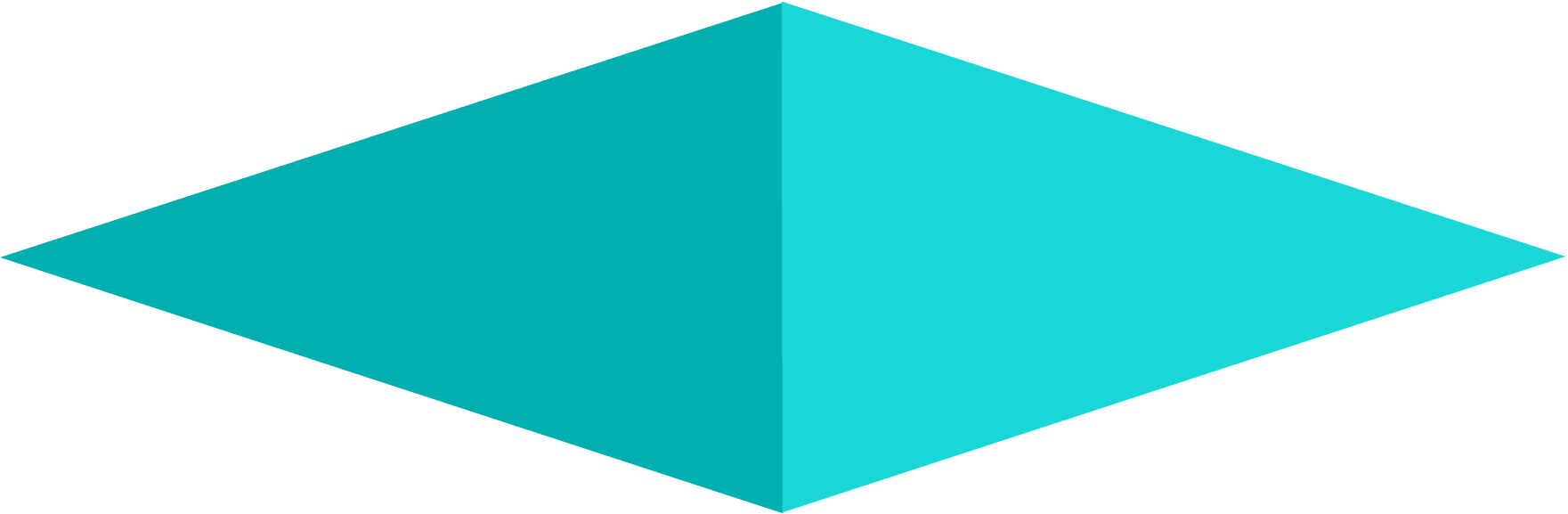}
\includegraphics[scale=0.3]{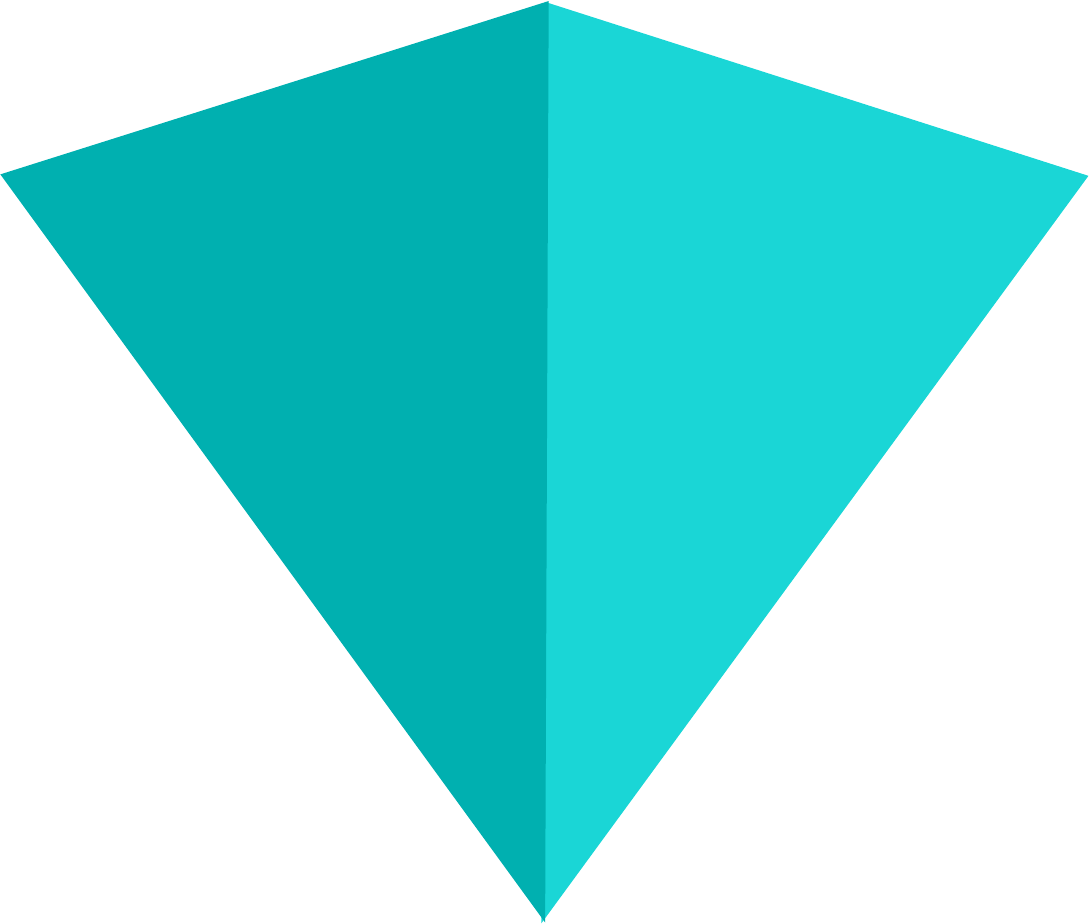}
\caption{Cutting Penrose tiles}
\label{halftiles}
\end{center}
\end{figure}

Decomposing Penrose tiles in half, yielding isosceles triangles as in Figure~\ref{halftiles}, is a very simple geometrical operation 
that has important repercussions. 

First of all, it is the first step of both the {\it inflation} 
and {\it deflation} procedures. In the case of inflation, the triangles are 
appropriately combined to form a new tiling, whose tiles are 
rescaled by a factor $\phi$. In the case of deflation, the triangles are further decomposed into smaller ones to yield the half--tiles 
of another tiling that is rescaled by a factor $1/\phi$. It is easy to see that these operations are
inverses of one another. We refer the reader to \cite{austin} for a detailed description in the case of rhombus tilings. 

Cutting kites in half can also be used to transform a kite and dart tiling into a rhombus tiling. 
\begin{figure}[h]
\begin{center}
\includegraphics[scale=0.8]{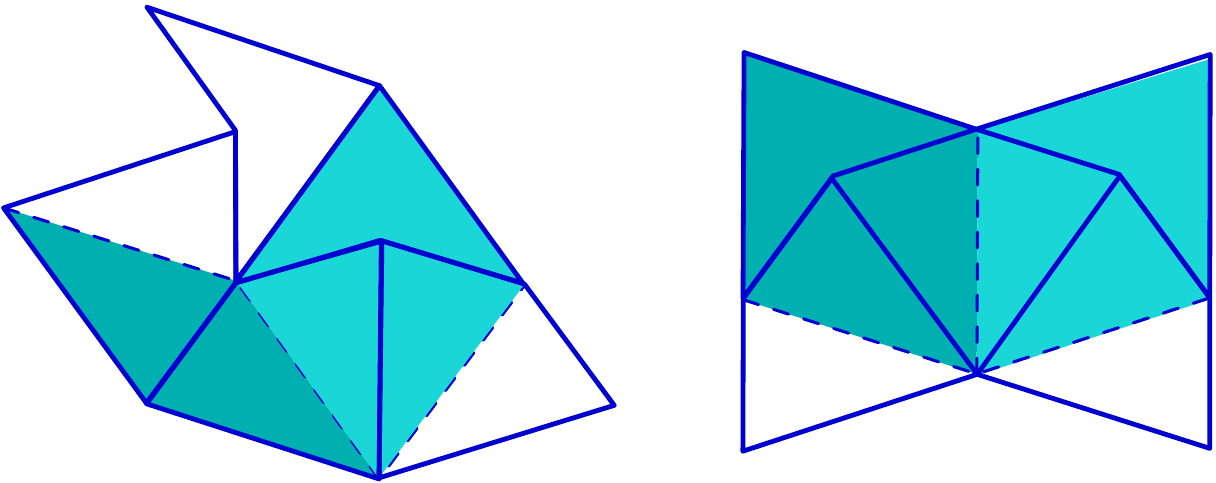}
\caption{From a kite and dart tiling to a rhombus tiling}
\label{passaggio}
\end{center}
\end{figure} 
The triangles are appropriately combined with each other and possibly a dart, in order to form thick and thin rhombuses
(see \cite{senechal1}  and Figure~\ref{passaggio}). 

Now, the process of subdividing a simple convex polytope into two smaller ones corresponds, 
at the (smooth) symplectic level, to the {\it symplectic cutting} operation, which was introduced by Lerman \cite{lerman}. 
In the toric setting, the original manifold decomposes into two new ones, 
each corresponding to one of the subdivided polytopes. 
The decomposition of Penrose tiles motivated us to extend this operation to the simple 
nonrational toric case. We find, for example, that the 
toric quasifold corresponding to each half--kite is given by
$$\frac{\left\{ \,(z_1,z_2,z_3)\in\C^3\,|\,|z_1|^2+\phi|z_2|^2+\phi |z_3|^2=1\,\right\}}{\left\{\,(e^{2\pi i s}, e^{2\pi i \phi s}, e^{2\pi i \phi(s+k)})\in \R^3/\Z^3\,|\,s\in\R,k\in \Z\,\right\}}.$$

\subsection*{The regular convex polyhedra}
The regular convex polyhedra are notable examples of convex polytopes and, as such, 
it is only natural to want to understand what the corresponding toric spaces look like. 
The cube and the regular tetrahedron are rational and simple, 
and yield smooth manifolds given, respectively, by $S^2\times S^2\times S^2$ and $\C\P^3$. 
The other three each present their complexities. The regular octahedron is rational but not simple, 
the regular dodecahedron is simple but not rational, while the regular icosahedron is neither rational nor simple. 
The first yields a space that is stratified by manifolds, the second yields a quasifold, 
while the third yields a space that is stratified by quasifolds; they are described explicitly in joint work with Battaglia.
The quasilattice for the dodecahedron is known in physics as the 
{\em simple icosahedral lattice} while the one for the icosahedron is known as the {\em body--centered icosahedral lattice}.
The normals, here too, are chosen among the quasilattice generators.

\subsection*{Final considerations}
\subsubsection*{Quasifolds, nonperiodic tilings and quasicrystals}
As we have shown, a number of interesting examples of toric quasifolds arise in connection with nonperiodic tilings.
There also appears to be a correspondence between some of the fundamental operations in the two theories.  
We have seen, in fact, that decomposing convex 
Penrose tiles into half corresponds to cutting the associated symplectic toric quasifolds. 
We expect, moreover, that recombining these half--tiles, as needed for the inflation and deflation procedures,
will correspond to a nonrational generalization of the inverse operation of symplectic cutting, which is given, 
in the smooth case, by the {\em symplectic sum} \cite{gompf}. We believe that it would be interesting 
to pursue the study of these connections even further. 
As a matter of fact, certain nonperiodic tilings have been used as 
mathematical models for the theory of {\em quasicrystals} \cite{senechal1}; these are special materials that were experimentally discovered by 
Shechtman et al. \cite{shechtman} that have discrete nonperiodic 
diffraction patterns. 
Actually, the very existence of these materials had been predicted by Mackay in connection with 
his studies of Penrose and Ammann tilings \cite{mackay1, mackay2}.
Ultimately, it is quite possible
that toric quasifolds might contribute to their theoretical understanding. A first step would consist in analyzing from the toric
viewpoint other tilings (and their operations) that are relevant in this respect.
Significant (though not the only) examples would be Socolar's octagonal and dodecagonal tilings, 
which are used as a basis for a treatment of the elasticity of octagonal and dodecagonal quasicrystals \cite{socolar}.  

\subsubsection*{Combinatorial equivalence in toric geometry}
By slightly perturbing the hyperplanes in (\ref{decomp}), it can be shown that every simple or simplicial polytope can 
be perturbed to a rational one that is combinatorially equivalent  (\cite[Proposition~2.17]{ziegler}). 
In the simple case, these perturbations yield, at the toric level, interesting families of quasifolds. For example, 
jointly with Battaglia and Zaffran, we have used one such perturbation to construct a one--parameter family of toric quasifolds 
that generalize and contain Hirzebruch surfaces.
This perturbation property also holds true for any three--dimensional polytope, 
not necessarily simple or simplicial \cite[Corollary~4.8]{ziegler}.
In greater dimensions, there are examples of polytopes for which this does not happen. 
The first was found by Perles in the sixties and has dimension $8$ (see \cite[Example~6.21]{ziegler}
and \cite{ziegler_int}).  As we have seen, from the toric viewpoint, these polytopes, being necessarily nonsimple, yield
spaces that are stratified by quasifolds. We believe it would be interesting to study these stratified 
spaces and understand how their geometry is affected by the fact that the corresponding polytopes 
cannot be deformed to rational ones within their combinatorial class.

\subsubsection*{Recent alternate approaches to nonrational toric geometry}
In recent years, there has been a surge of interest in nonrational toric geometry, 
and several alternate approaches to this subject have been introduced.
It should be said, first of all, that toric quasifolds can be thought of both as examples of stacks and of diffeological spaces. 
The stack approach to nonrational toric geometry was espoused first by Hoffman--Sjamaar
\cite{hoffmansjamaar, hoffman} and then by Katzarkov et al. \cite{KLMV}. 
Diffeological quasifolds, on the other hand, were studied jointly with
Iglesias--Zemmour in \cite{diffeoquasi}, 
providing an explicit link with non--commutative geometry \cite{connes}; applications of this viewpoint to the toric setting are work
in progress. Other recent points of view, due to Battaglia--Zaffran \cite{bz1, bz2}, Lin--Sjamaar \cite{linsjamaar}, Ratiu--Zung \cite{ratiuzung}, and Ishida et al. \cite{IKP}, involve foliations of smooth manifolds, either in the complex or presymplectic setting. 
We would like to point out that most of the
alternate nonrational toric viewpoints are founded on variations on the theme
of the fundamental triple (\ref{triple}), beginning, first and foremost, with the quasilattice $Q$. 
In joint work with Battaglia \cite{whatis}, we describe in detail how many of these different 
perspectives connect with each other and with ours; 
a dictionary is provided, in the hope that it will provide clarity and facilitate future interaction in the field.

\bigskip

\noindent 
Dipartimento di Matematica e Informatica "U. Dini", Universit\`a di Firenze \\
Piazza Ghiberti 27, 50122 Firenze, ITALY

\noindent
{\tt elisa.prato@unifi.it}
\end{document}